\documentclass{article}
\usepackage{amsmath,amsthm,amssymb}
\usepackage{url,fullpage}

\newtheorem{theorem}{Theorem}[section]
\theoremstyle{definition}
\newtheorem{remark}{Remark}[section]

\title{A Simple and Elementary Proof of Zorn's Lemma}
\author{Koji Nuida}
\date{Institute of Mathematics for Industry (IMI), Kyushu University \\ \url{nuida@imi.kyushu-u.ac.jp}}

\begin{document}

\maketitle

\begin{abstract}
Zorn's Lemma is a well-known equivalent of the Axiom of Choice.
It is usually regarded as a topic in axiomatic set theory, and its historically standard proof (from the Axiom of Choice) relies on transfinite recursion, a non-elementary set-theoretic machinery.
However, the statement of Zorn's Lemma itself uses only elementary terminology for partially ordered sets.
Therefore, it is worthy to establish a proof using only such elementary terminology.
Following this line of study, we give a new simple proof of Zorn's Lemma, which does not even use the notion of a well-ordered set.
\\[2mm]
\end{abstract}

\section{Introduction}

Zorn's Lemma is one of the most famous equivalents of the Axiom of Choice under the Zermelo--Fraenkel set theory.
It is historically regarded as a topic of axiomatic set theory, and the \lq\lq standard\rq\rq{} proof of Zorn's Lemma (from the Axiom of Choice) relies a non-elementary set-theoretic machinery called transfinite recursion (see e.g., \cite[Theorem 6.1]{Cameron}, \cite[Theorem 2.2]{Herrlich}, \cite[Theorem 5.4]{Jech_Set}, \cite[Theorem 2.1]{Jech_AC}).
However, Zorn's Lemma is ubiquitously used in mathematics, not just in set theory, and Zorn's Lemma can be stated with only elementary terminology for partially ordered sets (posets).
Therefore, it is worthy to establish a proof of Zorn's Lemma by using such elementary terminology only.
There exists a large number of proofs, without transfinite recursion, of Zorn's Lemma \cite{Kneser50,Szele50,Weston57,Halmos,Dugundji,Nagata,RubRub,Lewin91,Bergman97,Lang,Grayson07,Jain12,Tausk18,An19,Bukh21,IncTer24} or some property that easily implies Zorn's Lemma, e.g., Hausdorff Maximal Principle \cite{Frink52,Rudin,Kuttler,Brown10,Lektorean10,Spronk} and existence of fixed points for some mappings on a poset \cite{Bourbaki49,Berge,Brown70,Lang_functional,Zarouali-Darkaoui19,Shurman}.
Following this direction of study, in this short note we give yet another simple proof of Zorn's Lemma using elementary terminology only.

Among those previous proofs, our proof here is a descendant of a proof by Lewin \cite{Lewin91}.
Intuitively speaking, Lewin's proof constructed a \lq\lq highest\rq\rq{} chain in a given poset as the union of some family $\mathcal{C}$ of chains, which was defined by using the notion of a well-ordered set.
The nontrivial point of our work is to extract what is really required in the proof among the properties of well-ordered sets.
From this point of view, we introduce two properties named \lq\lq (i-$C$)\rq\rq{} and \lq\lq (ii-$C$)\rq\rq; a trick in our proof is that the construction of the family $\mathcal{C}$ is now in two steps, where the first step uses the condition (i-$C$) to construct an auxiliary family $\mathcal{C}_0$ and the second step constructs $\mathcal{C}$ by using the condition (ii-$C$) defined in terms of the family $\mathcal{C}_0$.
As a result, the terminology in our proof is more elementary, not even using the notion of a well-ordered set.

\paragraph{Notations and Terminology.}

In this note, a \emph{poset} (partially ordered set) is denoted by $(P,\leq)$ where $P$ is a set and $\leq$ is a binary relation on $P$ with the following three axioms (for arbitrary $x,y,z \in P$): (I) $x \leq x$; (II) $x \leq y$ and $y \leq x$ imply $x = y$; (III) $x \leq y$ and $y \leq z$ imply $x \leq z$.
We write $x < y$ to mean \lq\lq $x \leq y$ and $x \neq y$\rq\rq.
A subset $C$ of $P$ is called a \emph{chain} in $P$ if any pair $(x,y)$ of elements of $C$ satisfies either $x \leq y$ or $y \leq x$.
Note that any subset of a chain in $P$ is also a chain in $P$.
We say that $x \in P$ is an \emph{upper bound} of a chain $C \subseteq P$ if $y \leq x$ holds for any $y \in C$; and $x \in P$ is a \emph{strict upper bound} of $C$ if $y < x$ holds for any $y \in C$, or equivalently, $x$ is an upper bound of $C$ and $x \not\in C$.
A poset $(P,\leq)$ is said to be \emph{inductively ordered} if any chain in $P$ has an upper bound in $P$.
We say that an element $x \in P$ is \emph{maximal} if there exists no $y \in P$ satisfying $x < y$.
On the other hand, we say that an element $x$ of a chain $C$ in $P$ is the \emph{maximum element} of $C$, denoted by $\max C$, if $y \leq x$ holds for any $y \in C$ (i.e., $x$ is an upper bound of $C$).
With these terminology, Zorn's Lemma is stated as follows.

\begin{theorem}
[Zorn's Lemma]
Any inductively ordered poset $(P,\leq)$ has a maximal element.
\end{theorem}

\section{Our Proof}

We describe our proof of Zorn's Lemma.
Let $\mathcal{T} := \{ C \subseteq P : \mbox{$C$ is a chain in $P$} \}$.
For each $C \in \mathcal{T}$, define
\[
  \overline{U}_C := \{ x \in P : \mbox{$x$ is an upper bound of $C$} \} \,,\,
  U_C := \{ x \in P : \mbox{$x$ is a strict upper bound of $C$} \} = \overline{U}_C \setminus C \enspace.
\]
Then $\overline{U}_C = U_C \cup \{\max C\}$ if the maximum element $\max C$ of $C$ exists, and $\overline{U}_C = U_C$ otherwise.
Therefore,
\begin{equation}
  \label{eq:U_is_upwards_closed}
  \mbox{if $x \in \overline{U}_C$, $y \in P$, and $x < y$, then $y \in U_C$}
\end{equation}
(note that if $\max C < y$ then $y \not\in C$).
Moreover,
\begin{equation}
  \label{eq:lemma_for_U}
  \mbox{if $C_1,C_2 \in \mathcal{T}$ and $U_{C_1} \not\subseteq U_{C_2}$, then $C_1 \cap \overline{U}_{C_2} = \emptyset$}
\end{equation}
(indeed, if $x \in C_1 \cap \overline{U}_{C_2}$, then any $y \in U_{C_1}$ satisfies that $x < y$ by definition of $U_{C_1}$ and hence $y \in U_{C_2}$ by \eqref{eq:U_is_upwards_closed}, contradicting the assumption $U_{C_1} \not\subseteq U_{C_2}$).
Now the Axiom of Choice yields a choice function $f_0$ for the family of non-empty subsets $X$ of $P$; that is, $f_0(X) \in X$ for any such $X$.
Then we define a function $f \colon \mathcal{T} \to P$ satisfying that, for each $C \in \mathcal{T}$,
\[
  f(C) := f_0(U_C) \in U_C \subseteq P \setminus C \mbox{ if $U_C \neq \emptyset$} \,,\,
  f(C) := \max C \in C \cap \overline{U}_C \mbox{ if $U_C = \emptyset$}
\]
(note that when $U_C = \emptyset$, the assumption on $(P,\leq)$ being inductively ordered implies that $C$ has an upper bound $x \in P$ that is not strict, which satisfies $x \in C$ and must be the maximum element of $C$).
By the construction of $f$,
\begin{equation}
  \label{eq:choice_function}
  \mbox{if $C_1,C_2 \in \mathcal{T}$ and $U_{C_1} = U_{C_2} \neq \emptyset$, then $f(C_1) = f(C_2) \in U_{C_1}$ ($ = U_{C_2}$)} \enspace.
\end{equation}

Now let $\mathcal{C}_0$ denote the set of all $C \in \mathcal{T}$ satisfying the following condition:
\begin{quote}
  (i-$C$) $S \subseteq C$ and $U_S \not\subseteq U_C$ imply $f(S) \in C$.
\end{quote}
Then let $\mathcal{C}$ denote the set of all $C \in \mathcal{C}_0$ satisfying the following condition:
\begin{quote}
  (ii-$C$) $C' \in \mathcal{C}_0$ implies $C \subseteq C' \cup U_{C'}$.
\end{quote}

Let $C^* := \bigcup \mathcal{C}$.
We prove that $C^* \in \mathcal{C}$, by verifying the defining conditions of $\mathcal{C}$ as follows:
\begin{itemize}
  \item 
  For (ii-$C^*$), let $C' \in \mathcal{C}_0$.
  Then each $C \in \mathcal{C}$ satisfies $C \subseteq C' \cup U_{C'}$ by (ii-$C$), therefore $C^* = \bigcup \mathcal{C} \subseteq C' \cup U_{C'}$, as desired.
  \item 
  We show that $C^* \in \mathcal{T}$, i.e., any $x,y \in C^*$ satisfy $x \leq y$ or $y \leq x$.
  We can take $C,C' \in \mathcal{C}$ with $x \in C$ and $y \in C'$.
  Now $C \subseteq C' \cup U_{C'}$ by (ii-$C$), therefore we have either $x,y \in C'$ or $x \in U_{C'}$ (hence $y < x$), implying the claim in any case.
  \item 
  For (i-$C^*$), suppose that $S \subseteq C^*$ and $U_S \not\subseteq U_{C^*} = \bigcap_{C \in \mathcal{C}} U_C$.
  Then $U_S \not\subseteq U_C$ for some $C \in \mathcal{C}$.
  Now $S \cap U_C = \emptyset$ by \eqref{eq:lemma_for_U}, while $S \subseteq C^* \subseteq C \cup U_C$ by (ii-$C^*$) applied to $C \in \mathcal{C}_0$, therefore $S \subseteq C$.
  By (i-$C$) applied to $S \subseteq C$, we have $f(S) \in C \subseteq C^*$, therefore $f(S) \in C^*$, as desired.
\end{itemize}

Now if $U_{C^*} = \emptyset$, then $f(C^*) = \max C^*$ is a maximal element of $P$ (as otherwise $\max C^* < y$ for some $y \in P$ and hence $y \in U_{C^*}$ by \eqref{eq:U_is_upwards_closed}, a contradiction), as desired.
Hence the proof will be completed once we obtain a contradiction assuming that $U_{C^*} \neq \emptyset$ (hence $f(C^*) \in U_{C^*}$).
Let $u := f(C^*) \in U_{C^*}$ and $C^{**} := C^* \cup \{u\}$; therefore we have $u = \max C^{**} \not\in C^*$, $C^{**} \not\subseteq C^*$, and $C^{**} \in \mathcal{T}$ as $C^* \in \mathcal{T}$.
We prove that $C^{**} \in \mathcal{C}$, by verifying the remaining defining conditions as follows:
\begin{itemize}
  \item
  For (i-$C^{**}$), suppose that $S \subseteq C^{**}$ and $U_S \not\subseteq U_{C^{**}}$.
  Then $S \cap \overline{U}_{C^{**}} = \emptyset$ by \eqref{eq:lemma_for_U} and hence $u \not\in S$.
  This implies that $S \subseteq C^*$ and hence $U_{C^*} \subseteq U_S$.
  Now if $U_S \subseteq U_{C^*}$, then $U_S = U_{C^*} \neq \emptyset$ and $f(S) = f(C^*) = u \in C^{**}$ by \eqref{eq:choice_function}; while if $U_S \not\subseteq U_{C^*}$, then we have $f(S) \in C^* \subseteq C^{**}$ by applying (i-$C^*$) to $S \subseteq C^*$.
  Hence $f(S) \in C^{**}$ in any case, as desired.
  \item
  For (ii-$C^{**}$), let $C' \in \mathcal{C}_0$.
  As $C^* \subseteq C' \cup U_{C'}$ by (ii-$C^*$), it suffices to show that $u \in C' \cup U_{C'}$, or equivalently, $u \in C'$ if $u \not\in U_{C'}$.
  Now we have $U_{C^*} \not\subseteq U_{C'}$ as $u \in U_{C^*}$.
  Hence $C^* \cap U_{C'} = \emptyset$ by \eqref{eq:lemma_for_U}, while $C^* \subseteq C' \cup U_{C'}$ as above, therefore $C^* \subseteq C'$.
  Hence $u = f(C^*) \in C'$ by applying (i-$C'$) to $C^* \subseteq C'$, as desired.
\end{itemize}
However, this fact $C^{**} \in \mathcal{C}$ yields a contradiction, as $C^{**} \not\subseteq C^* = \bigcup \mathcal{C}$.
This completes the proof.
\qed

\begin{remark}
  If we change the family $\mathcal{T}$ in our proof to the family of all well-ordered subsets of $P$, then it yields a proof with weakened assumption on $(P,\leq)$ being inductively ordered where the existence of upper bounds is now assured only for well-ordered subsets of $P$.
  Indeed, now for showing that the set $C^*$ in the proof is a member of $\mathcal{T}$, for any non-empty subset $S \subseteq C^*$, fix any $C \in \mathcal{T}$ with $S \cap C \neq \emptyset$, and take $x := \min (S \cap C)$.
  For $y \in S$, if $y \in C$ then $x = \min (S \cap C) \leq y$; while if $y \not\in C$ then $y \in C^* \setminus C \subseteq U_C$ from \emph{(ii-$C^*$)} and $x < y$.
  Hence $x = \min S$, therefore $C^* \in \mathcal{T}$.
  The proof of $C^{**} \in \mathcal{T}$ is almost the same as the original proof, and the remaining part of the proof is not affected by the change of the definition of $\mathcal{T}$.
\end{remark}

We explain the difference of our proof from Lewin's work \cite{Lewin91} mentioned in the Introduction.
The outline of the proof is common to both proofs, i.e., (1) defining some family $\mathcal{C}$ of chains and (2) showing that $C^* := \bigcup \mathcal{C} \in \mathcal{C}$ and that $C^{**} := C^* \cup \{f(C^*)\} \in \mathcal{C}$ where $f(C^*)$ is a strict upper bound of $C^*$, yielding a contradiction.
However, in contrast to our proof where the key property $C^* \in \mathcal{C}$ in Step (2) is derived directly from the defining conditions for $\mathcal{C}$, Lewin's proof required an intermediate step to show some extra property for $\mathcal{C}$ that is seemingly stronger than the defining conditions for $\mathcal{C}$.
In detail, the following comparability property was shown in Lewin's proof: for any two members of $\mathcal{C}$, one of them is an initial segment of the other.
This is analogous to a property of well-ordered sets, and to ensure this property, $\mathcal{C}$ was defined in a way that each member of $\mathcal{C}$ should be a well-ordered subset of $P$.
Our main idea is that the full comparability property is in fact not necessary in the proof.
Our new condition (ii-$C$) can be seen as a weaker variant of the comparability property (indeed, the condition $C \subseteq C' \cup U_{C'}$ holds when one of $C$ and $C'$ is an initial segment of the other), which (for $C^* = \bigcup \mathcal{C}$ and $C^{**} = C^* \cup \{f(C^*)\}$) can be proved directly without requiring that members of $\mathcal{C}$ are well-ordered subsets.

\paragraph{Acknowledgment.}

The author deeply thanks \'{A}rp\'{a}d Sz\'{a}z for several detailed comments to a previous version of the manuscript, including a comprehensive list of known proofs of Zorn's Lemma without transfinite recursion (on which the reference list in the current paper is based).
The author also thanks Pedro S\'{a}nchez Terraf for kindly informing the author of another reference \cite{IncTer24} for a proof of Zorn's Lemma.

%%%%%%%%%%%%%%%%%%%%%%%%%%%%%%%%%%%%%%%%%%%%%%%%%%%%%%%%%%%%%%%%%%%%

\end{document}